\documentclass[runningheads,a4paper]{llncs}
\usepackage{a4wide}
\usepackage{algorithmic}
\usepackage{amsfonts,amsbsy}
\usepackage{amssymb}
\usepackage{amstext}
\usepackage{color}

\def\proof{\par\noindent{\bf Proof:~}}
\def\QED{\quad\hbox{\hskip 1pt \vrule width 4pt height 8pt depth 1.5pt
         \hskip 1pt}\lower 8.5pt\null\par}
\def\inline#1:{\par\vskip 7pt\noindent{\bf #1:}\hskip 10pt}
\newcommand{\bfm}{\textbf}
\newcommand{\mcal}{\mathcal}
\newcommand{\mbb}{\mathbb}
\newcommand{\mbs}{\boldsymbol}

\newcommand{\E}{\mcal{E}}
\newcommand{\R}{\mcal{R}}
\newcommand{\A}{\mcal{A}}

\newcommand{\prob}[1]{\ensuremath{\text{{\bf Pr}$\left[#1\right]$}}}
\newcommand{\expect}[1]{\ensuremath{\text{{\bf E}$\left[#1\right]$}}}

\newtheorem{observation}[theorem]{Observation}
\newcommand{\etal}{\textit{et~al.}}

\begin{document}
\mainmatter

\title{Conflict-Free Coloring Made Stronger}
\titlerunning{Conflict-Free Coloring Made Stronger}

\author{
Elad Horev\thanks{Computer Science department, Ben-Gurion
University, Beer Sheva, Israel; {\tt horevel@cs.bgu.ac.il}.} \and
Roi Krakovski\thanks{Computer Science department, Ben-Gurion
University, Beer Sheva, Israel; {\tt roikr@cs.bgu.ac.il}.} \and
Shakhar Smorodinsky\thanks{Mathematics department, Ben-Gurion
University, Beer Sheva, Israel; {\tt
http://www.math.bgu.ac.il/$\sim$shakhar/ ;
shakhar@math.bgu.ac.il}.}}

\authorrunning{Horev, Krakovski, and Smorodinsky}
\toctitle{Lecture Notes in Computer Science}
%\tocauthor{}

\institute{Ben-Gurion University, Beer Sheva, Israel}
\maketitle

\begin{abstract}
In FOCS 2002, Even et al. showed that any set of $n$ discs in the
plane can be Conflict-Free colored with a total of at most
$O(\log n)$ colors. That is, it can be colored with $O(\log n)$
colors such that for any (covered) point $p$ there is some disc
whose color is distinct from all other colors of discs containing
$p$. They also showed that this bound is asymptotically tight. In
this paper we prove the following stronger results:

\begin{enumerate}
\item [(i)] Any set of $n$ discs in the plane can be colored with a total
of at most $O(k \log n)$ colors such that (a) for any point $p$ that
is covered by at least $k$ discs, there are at least $k$ distinct
discs each of which is colored by a color distinct from all other
discs containing $p$ and (b) for any point $p$ covered by at most $k$
discs, all discs covering $p$ are colored distinctively. We call
such a coloring a {\em $k$-Strong Conflict-Free} coloring. We
extend this result to pseudo-discs and arbitrary regions with
linear union-complexity.

\item [(ii)] More generally, for families of $n$ simple closed Jordan
regions with union-complexity bounded by $O(n^{1+\alpha})$, we
prove that there exists a $k$-Strong Conflict-Free coloring with
at most $O(k n^\alpha)$ colors.

\item [(iii)] We prove that any set of $n$ axis-parallel rectangles can be
$k$-Strong Conflict-Free colored with at most $O(k \log^2 n)$
colors.

\item [(iv)] We provide a general framework for $k$-Strong Conflict-Free
coloring arbitrary hypergraphs. This framework relates the notion
of $k$-Strong Conflict-Free coloring and the recently studied
notion of $k$-colorful coloring.
\end{enumerate}

All of our proofs are constructive. That is, there exist
polynomial time algorithms for computing such colorings.\\

\noindent
\bfm{Key Words.} Conflict-Free Colorings, Geometric hypergraphs, Wireless networks, Discrete geometry. 
\end{abstract}

\section{Introduction and Preliminaries}
\label{sec:intro} Motivated by modeling frequency assignment to
cellular antennae, Even et al. \cite{even} introduced the notion
of Conflict-Free colorings. A {\em Conflict-Free} coloring (CF in
short) of a hypergraph $H=(V,\E)$ is a coloring of the vertices
$V$ such that for any non-empty hyperedge $e \in \E$ there is
some vertex $v \in e$ whose color is distinct from all other
colors of vertices in $e$. For a hypergraph $H$, one seeks the
least number of colors $l$ such that there exists an $l$-coloring
of $H$ which is Conflict-Free. It is easily seen that CF-coloring
of a hypergraph $H$ coincides with the notion of classical graph coloring in
the case when $H$ is a graph (i.e., all hyperedges are of
cardinality two). Thus it can be viewed as a
generalization of graph coloring. There are two well known
generalizations of graph coloring to hypergraph coloring in the
literature (see, e.g., \cite{BERGE}). The first generalization
requires ``less" than the CF requirement and this is the {\em
non-monochromatic} requirement where each hyperedge in $\E$ of
cardinality at least two should be non-monochromatic: The
\emph{chromatic number} of a hypergraph $H$, denoted $\chi(H)$,
is the least number $l$ such that $H$ admits an $l$-coloring
which is a non-monochromatic coloring. The second generalization
requires ``more" than the CF requirement and this is the {\em
colorful} requirement where each hyperedge should be colorful
(i.e., all of its vertices should have distinct colors). For
instance, consider the following hypergraph $H=(V,\E)$: Let $V =
\{1,2,\ldots,n\}$ and let $\E$ consist of all subsets of $V$
consisting of consecutive numbers of $V$. That is, $\E$ consists of
all discrete intervals of $V$. It is easily seen that one can
color the elements of $V$ with two colors in order to obtain a
non-monochromatic coloring of $H$. Color the elements of $V$
alternately with `black' and `white'. On the other extreme, one
needs $n$ colors in any colorful coloring of $H$. Indeed $V$
itself is also a hyperedge in this hypergraph (an `interval'
containing all elements of $V$) so all colors must be distinct.
However, it is an easy exercise to see that there exists a
CF-coloring of $H$ with $\lfloor \log n \rfloor +1$ colors. In
fact, for an integer $k>0$, if $V$ consist of $2^k -1$ elements
then $k$ colors suffice and are necessary for CF-coloring $H$.

Let $\R$ be a finite collection of regions in $\mbb{R}^d$, $d
\geq 1$. For a point $p \in \mbb{R}^d$, define $r(p) = \{R \in
\R: p \in R\}$. The hypergraph $(\R, \{r(p)\}_{p\in \mbb{R}^d})$,
denoted $H(\R)$, is called the hypergraph \emph{induced} by $\R$.
Such hypergraphs are referred to as \emph{geometrically induced}
hypergraphs. Informally these are the Venn diagrams of the
underlying regions.

In general, dealing with CF coloring for arbitrary hypergraphs is
not easier than graph coloring. The paper \cite{even} focused on
hypergraphs that are induced by geometric objects such as discs,
squares etc. Their motivation was a modeling of frequency
assignment to cellular antennae in a manner that reduces the
spectrum of frequencies used by a network of antennae. Suppose
that antennae are represented by discs in the plane and that every
client (holding a cell-phone) is represented by a point. Antennae
are assigned frequencies (this is the coloring). A client is
served provided that there is at least one antenna `covering' the
client for which the assigned frequency is ``unique" and
therefore has no ``conflict" (interference) with other
frequencies used by nearby antennae. When $\R$ is a finite family
of $n$ discs in the plane $\mbb{R}^2$, Even et al. \cite{even}
proved that finding an optimal CF-coloring for $\R$ is NP-hard.
However, they showed that there is always a CF-coloring of
$H(\R)$ with $O(\log n)$ colors and that this bound is
asymptotically tight. That is, for every $n$ there is a family of
$n$ discs which requires $\Omega(\log n)$ colors in any
CF-coloring. See \cite{even} for further discussion of this model
and the motivation.

CF-coloring finds application also in activation protocols for
RFID networks. Radio frequency identification (RFID) is a
technology where a reader device can "sense" the presence of a
nearby object by reading a tag device attached to the object. To
improve coverage, multiple RFID readers can be deployed in the
given region. However, two readers trying to access a tagged
device simultaneously might cause mutual interference. One may
want to design scheduled access of RFID tags in a multiple reader
environment. Assume that we have $t$ time slots and we would like
to `color' each reader with a time slot in $\{1,\ldots,t\}$ such
that the reader will try to read all nearby tags in its given
time slot. In particular, we would like to read all the tags and
minimize the total time slots $t$. It is easily seen that if we
CF-color the family $\R$ of readers then in this coloring every
possible tag will have a time slot and a single reader trying to
access it in that time slot \cite{HIMANSHU}. The notion of
CF-coloring has caught much scientific attention in recent years
both from the algorithmic and combinatorial point of view
\cite{cf9,cf8,BarNoy3,BarNoy2,NCS06,cf7,cf6,cf5,cf4,cf3,cf2,PT09,cf1,smoro}.

\paragraph{Our Contribution:}
 In this paper we study the notion of
\emph{$k$-Strong-Conflict-Free} (abbreviated, $kSCF$) colorings
of hypergraphs. This notion extends the  notion of $CF$-colorings
of hypergraphs. Informally, in the case of coloring discs, rather
than having at least one unique color at every covered point $p$,
we require at least $k$ distinct colors to some $k$ discs such
that each of these colors is unique among the discs covering $p$.
The motivation for studying $kSCF$-coloring is rather
straightforward in the context of wireless antennae. Having, say
$k>1$ unique frequencies in any given location allows us to serve
$k$ clients at that location rather than only one client. In the
context of RFID networks, a $k$SCF coloring will correspond to an
activation protocol which is fault-tolerant. That is, every tag
can be read even if some $k-1$ readers are broken.

\begin{definition}[$\mbs{k}$-Strong Conflict-Free coloring:]
Let $H= (V, \E)$ be a hypergraph and let  $k \in \mbb{N}$ be some
fixed integer. A coloring of $V$ is called {\em
$k$-Strong-Conflict-Free} for $H$ if\\ 
(i) for every hyperedge $e \in \E$ with $|e| \geq k$ there exists at least $k$ vertices in $e$,
whose colors are unique among the colors assigned to the vertices of $e$, and\\
(ii) for each hyperedge $e \in \E$ with $|e| < k$ all vertices in $e$ get distinct colors.

Let $f_H(k)$ denote the least integer $l$ such that $H$ admits a $kSCF$-coloring with $l$
colors.
\end{definition}
Note that a $CF$-coloring of a hypergraph $H$ is $kSCF$-coloring
of $H$ for $k=1$.

Abellanas et al. \cite{HURTADO} were the first to study
$k$SCF-coloring\footnote{They referred to such a coloring as
$k$-Conflict-Free coloring.}. They focused on the special case
where $V$ is a finite set of points in the plane and $\E$ consist
of all subsets of $V$ which can be realized as an intersection of
$V$ with a disc. They showed that in this case the hypergraph
admits a $k$SCF-coloring with $O(\frac{\log n}{\log
{\frac{ck}{ck-1}}})$ ($ = O(k \log n)$) colors, for some absolute
constant $c$. See also \cite{Abam08}.\\

The following notion of \emph{$k$-colorful} colorings was recently introduced and studied by
Aloupis et al. \cite{aloupis} for the special case of hypergraphs induced by discs.

\begin{definition}
\label{def:colorful} Let $H= (V, \E)$ be a hypergraph, and let
$\varphi$ be a coloring of $H$. A hyperedge $e \in \E$ is said to
be \emph{$k$-colorful} with respect to $\varphi$ if there exist
$k$ vertices in $e$ that are colored distinctively under
$\varphi$. The coloring $\varphi$ is called {\em $k$-colorful} if
every hyperedge $e \in \E$ is $\min \{|e|,k\}$-colorful. Let
$c_H(k)$ denote the least integer $l$ such that $H$ admits a
$k$-colorful coloring with $l$ colors.
\end{definition}

Aloupis et al. were motivated by a problem related to battery
lifetime in sensor networks. See \cite{aloupis,efrat,PachToth}
for additional details on the motivation and related problems.\\

{\bf \noindent Remark:} Every $kSCF$-coloring of a hypergraph $H$
is a $k$-colorful coloring of $H$. However, the opposite claim is
not necessarily true. A $k$-colorful coloring assures us that
every hyperedge of cardinality at least $k$ has at least $k$
distinct colors present in it. However, these $k$ colors are not
necessarily unique since each may appear with multiplicity.\\

A $k$-colorful coloring can be viewed as a type of coloring
which is ``in between" non-monochromatic coloring and colorful
coloring. A $2$-colorful coloring of $H$ is exactly the classical
non-monochromatic coloring, so $\chi(H) = c_H(2)$. If $H$ is a
hypergraph with $n$ vertices, then an $n$-colorful coloring of
$H$ is the classical colorful coloring of $H$.  Consider the
hypergraph $H$, consisting of all discrete intervals on
$V=\{1,\ldots,n\}$ mentioned earlier. It is easily seen that for
any $i$, an $i$-colorful coloring with $i$ colors is obtained by
coloring $V$ in increasing order with
$1,2,\ldots,i,1,2,\ldots,i,1,2\ldots$ with repetition.

In this paper, we study a connection between $k$-colorful coloring
and Strong-Conflict-Free coloring of hypergraphs. We show that if
a hypergraph $H$ admits a $k$-colorful coloring with a ``small"
number of colors (hereditarily) then it also admits a
$(k-1)$SCF-coloring with a ``small" number of colors. The
interrelation between the quoted terms is provided in
Theorems~\ref{th:reduce} and \ref{th:nonlinear} below.

Let $H=(V,\E)$ be a hypergraph and let $V' \subset V$. We write
$H[V']$ to denote the sub-hypergraph of $H$ induced by $V'$, i.e.,
$H[V'] = (V',\E')$ and $\E' = \{e \cap V' |e \in \E\}$. We write
$n(H)$ to denote the number of vertices of $H$.

\begin{theorem}
\label{th:reduce} Let $H = (V, \E)$ be a hypergraph with $n$
vertices, and let $k,\ell \in \mbb{N}$ be fixed integers, $k \geq
2$. If every induced sub-hypergraph $H'\subseteq H$ satisfies
$c_{H'}(k) \leq \ell$, then $f_H(k-1) \leq
\log_{1+\frac{1}{\ell-1}} n = O(l \log n)$.
\end{theorem}

\begin{theorem}
\label{th:nonlinear} Let $H = (V, \E)$ be a hypergraph with $n$
vertices, let $k \geq 2$ be a fixed integer, and let $0 < \alpha
\leq 1$ be a fixed real. If every induced sub-hypergraph
$H'\subseteq H$ satisfies $c_{H'}(k) = O(k{n(H')}^{\alpha})$,
then $f_H(k-1) = O(k{n(H')}^{\alpha})$.
\end{theorem}

Consider the hypergraph of ``discrete intervals" with $n$
vertices. As mentioned earlier, it has a $(k+1)$-colorful
coloring with $k+1$ colors and this holds for every induced
sub-hypergraph. Thus, Theorem~\ref{th:reduce} implies that it
also admits a $kSCF$-coloring with at most $\log_{1+\frac{1}{k}} n
= O(k \log n)$ colors. In Section~\ref{sec:upper}, we provide an
upper bound on the number of colors required by $kSCF$-coloring
of geometrically induced hypergraphs as a function of the
union-complexity of the regions that induce the hypergraphs.
Below we describe the relations between the union-complexity of
the regions, $k$-colorful and $(k-1)$SCF coloring of the
underlying hypergraph. First, we need to define the notion of
union-complexity.

\begin{definition}
For a family $\R$ of $n$ simple closed Jordan regions in the
plane, let $\partial{\R}$ denote the boundary of the union of the
regions in $\R$. The {\em union-complexity} of $\R$ is the number
of intersection points, of a pair of boundaries of regions in
$\R$, that belong to $\partial{\R}$.
\end{definition}

For a set $\R$ of $n$ simple closed planar Jordan regions, let
$\mcal{U}_{\R}:\mbb{N} \rightarrow \mbb{N}$ be a function such
that $\mcal{U}_{\R}(m)$ is the maximum union-complexity of any
subset of $k$ regions in $\R$ over all $k \leq m$, for $1 \leq m
\leq n$. We abuse the definition slightly and assume that the
union-complexity of any set of $n$ regions is at least $n$. When
dealing with geometrically induced hypergraphs, we consider
$k$-colorful coloring and $kSCF$-coloring of hypergraphs that are
induced by simple closed Jordan regions having union-complexity at
most $O(n^{1+\alpha})$, for some fixed parameter $0 \leq \alpha
\leq 1$. The value $\alpha =0$ corresponds to regions with linear
union-complexity such as discs or pseudo-discs (see, e.g.,
\cite{KLPS}). The value $\alpha =1$ corresponds to regions with
quadratic union-complexity. See \cite{Efrat2,Efrat1} for
additional families with sub-quadratic union-complexity.

In the following theorem we provide an upper bound on the number
of colors required by a $k$-colorful coloring of a geometrically
induced hypergraph as a function of $k$ and of the
union-complexity of the underlying regions inducing the
hypergraph :

\begin{theorem}
\label{th:prior-res} Let $k \geq 2$, let $0 \leq \alpha \leq 1$,
and let $c$ be a fixed constant. Let $\R$ be a set of $n$ simple
closed Jordan regions such that $\mcal{U}_{\R}(m) \leq c
m^{1+\alpha}$, for $1 \leq m \leq n$, and let $H=H(\R)$. Then
$c_H(k) = O( k n^{\alpha} )$.
\end{theorem}

Combining Theorem~\ref{th:reduce} with Theorem~\ref{th:prior-res}
(for $\alpha = 0$) and Theorem~\ref{th:nonlinear} with
Theorem~\ref{th:prior-res} (for $0 < \alpha < 1$) yields the
following result:

\begin{theorem}
\label{th:results} Let $k \geq 2$, let $0 \leq \alpha \leq 1$,
and let $c$ be a constant. Let $\R$ be a set of $n$ simple closed
Jordan regions such that $\mcal{U}_{\R}(m)= c m^{1+\alpha}$, for
$1 \leq m \leq n$. Let $H=H(\R)$. Then:
$$
 f_H(k-1) = \left\{ \begin{array}{ll}
                   O(k \log n), \mbox{      $\alpha = 0$}, \\
                   O(kn^{\alpha}), \mbox{       $0 < \alpha \leq 1$}.
                   \end{array}
            \right.
$$
\end{theorem}

In Section~\ref{sec:rects} we consider $kSCF$-colorings of
hypergraphs induced by axis-parallel rectangles in the plane. It
is easy to see that axis-parallel rectangles might have quadratic
union-complexity, for example, by considering a grid-like
construction of $n/2$ disjoint (horizontally narrow) rectangles
and $n/2$ disjoint (vertically narrow) rectangles. For a
hypergraph $H$ induced by  axis-parallel rectangles,
Theorem~\ref{th:results} states that $f_H(k-1) =O(k n)$. This
bound is meaningless, since the bound $f_H(k-1) \leq n$ is
trivial. Nevertheless, we provide a near-optimal upper bound for
this case in the following theorem:

\begin{theorem}
\label{th:rects} Let $k \geq 2$. Let $\R$ be a set of $n$
axis-parallel rectangles, and let $H=H(\R)$. Then $f_H(k-1) =O(k
\log^2 n)$.
\end{theorem}

In order to obtain Theorem~\ref{th:rects} we prove the following
theorem:

\begin{theorem} \label{th:colorful_rects}
Let $H = H(\R)$ be the
hypergraph induced by a family $\R$ of $n$ axis-parallel
rectangles in the plane, and let $k \in \mbb{N}$ be an integer,
$k \geq 2$. For every induced sub-hypergraph $H'\subseteq H$ we
have: $c_{H'}(k) \leq k \log n$.
\end{theorem}

Theorem~\ref{th:rects} is therefore an easy corollary
of Theorem~\ref{th:colorful_rects} combined with
Theorem~\ref{th:reduce}. 

Har-Peled and Smorodinsky \cite{cf3} proved that any
family $\R$ of $n$ axis-parallel rectangles admits a CF-coloring
with $O(\log ^2 n)$ colors. Their proof uses the probabilistic
method. They also provide a randomized algorithm for obtaining
CF-coloring with at most $O(\log ^2 n)$ colors. Later,
Smorodinsky \cite{smoro} provided a deterministic polynomial-time
algorithm that produces a CF-coloring for $n$ axis-parallel
rectangles with $O(\log^2 n)$ colors. Theorem~\ref{th:rects} thus
generalizes the results of \cite{cf3} and \cite{smoro}.

All of our proofs are constructive. In other words, there exist
deterministic polynomial-time algorithms to obtain the required
$kSCF$ coloring with the promised bounds. In this paper, we omit
the technical details of the underlying algorithms and we do not
make an effort to optimize their running time.

The result of Ali-Abam \etal \cite{Abam08} implies that the upper
bounds provided in Theorem~\ref{th:results} for $\alpha = 0$ and
Theorem~\ref{th:rects} are optimal. Specifically, they provide
matching lower bounds on the number of colors required by any
$kSCF$-coloring of hypergraphs induced by (unit) discs and
axis-parallel squares in the plane by a simple analysis of such
coloring for the discrete intervals hypergraph mentioned earlier.\\

\noindent
\bfm{Organization.} In Section~\ref{sec:general} we prove Theorems~\ref{th:reduce} and \ref{th:nonlinear}. In Section~\ref{sec:upper} we prove Theorems~\ref{th:prior-res} and \ref{th:results}. Finally, in Section~\ref{sec:rects} we prove Theorems~\ref{th:rects} and \ref{th:colorful_rects}.

\section{A Framework For Strong-Conflict-Free Coloring}
\label{sec:general}

In this section, we prove Theorems~\ref{th:reduce} and \ref{th:nonlinear}. 
To that end we devise a framework for obtaining an upper bound on the number
of colors required by a Strong-Conflict-Free coloring of a
hypergraph. Specifically, we show that if there
exist fixed integers $k$ and $l$ such that an $n$-vertex
hypergraph $H$ admits the hereditary property that every
vertex-induced sub-hypergraph $H'$ of $H$ admits a $k$-colorful
coloring with at most $l$ colors, then $H$ admits a
$(k-1)SCF$-coloring with $O(l \log n)$ colors. For the case when
$l$ is replaced with the function $k{n(H')}^{\alpha}$ we get a
better bound without the $\log n$ factor.

\noindent
\bfm{Framework $\mbs{\A}$:}\\
\bfm{Input:} A hypergraph $H$ satisfying the conditions of Theorems~\ref{th:reduce} and \ref{th:nonlinear}.\\
\bfm{Output:} A $(k-1)SCF$-coloring of $H$.
\begin{algorithmic}[1]
    \STATE $i \leftarrow 1$ \COMMENT{$i$ denotes an unused color.}
    \WHILE{$V \not= \emptyset$}
        \STATE \bfm{Auxiliary Coloring:} Let $\varphi:V \rightarrow [\ell]$ be a $k$-colorful coloring of $H[V]$ with at most $\ell$ colors.
        \STATE Let $V'$ be a color class of $\varphi$ of maximum cardinality.
        \STATE \bfm{Color:} Set $\chi(u)=i$ for every vertex $u \in V'$.
        \STATE \bfm{Discard:} $V \leftarrow V \setminus V'$.
        \STATE \bfm{Increment:} $i \leftarrow i+1$.
    \ENDWHILE
    \STATE \bfm{Return} $\chi$.
\end{algorithmic}
\ \\
\noindent
\bfm{Proof of Theorems~\ref{th:reduce} and~\ref{th:nonlinear}.}
We show that the coloring produced by Framework $\A$ is a $(k-1)$SCF-coloring of $H$ with a total number of colors as specified in Theorems~\ref{th:reduce} and~\ref{th:nonlinear}.

Let $\chi$ denote the coloring obtained by the application of
framework $\A$ on $H$. The number of colors used by $\chi$ is the
number of iterations performed by $\A$. By the pigeon-hole
principle, at least $|V|/\ell$ vertices are removed in each
iteration (where $V$ is the set of vertices remained after the
last iteration). Therefore, the total number of iterations
performed by $\A$ is bounded by $\log_{1+\frac{1}{\ell-1}} n$.
Thus, the coloring $\chi$ uses at most $\log_{1+\frac{1}{\ell-1}}
n$ colors. If in step 3 of the framework, $l$ is replaced with the
function $k{|V|}^{\alpha}$ (for a fixed parameter $0 < \alpha <
1$), then by the pigeon-hole principle at least
$\frac{|V|^{1-\alpha}}{k}$ vertices of $H$ are discarded in step
6 of that iteration. It is easily seen that the number of
iterations performed in this case is bounded by $O(kn^{\alpha})$
where $n=n(H)$.

Next, we prove that the coloring $\chi$ is indeed a
$(k-1)SCF$-coloring of $H$. The colors of $\chi$ are the indices
of iterations of $\A$. Let $e \in \E$ be a hyperedge of $H$. If
$|e| \leq k$ then it is easily seen that all colors of vertices
of $e$ are distinct. Indeed, by the property of the auxiliary
coloring $\varphi$ in step $3$ of the framework, every vertex of
$e$ is colored distinctively and in each such iteration, at most
one vertex from $e$ is colored by $\chi$ so $\chi$ colors all
vertices of $e$ in distinct iterations. Next, assume that $|e| >
k$. We prove that $e$ contains at least $k-1$ vertices that are
assigned unique colors in $\chi$. For an integer $r$, let
$\{\alpha_1,\ldots,\alpha_r\}$ denote the $r$ largest colors in
decreasing order that are assigned to some vertices of $e$. That
is, the color $\alpha_1$ is the largest color assigned to a
vertex of $e$, the color $\alpha_2$ is the second largest color
and so on. In what follows, we prove a stronger assertion that
for every $1 \leq j \leq k-1$ the color $\alpha_j$ exists and is
unique in $e$. The proof is by induction on $j$. $\alpha_1$
exists in $e$ by definition.  For the base of the induction we
prove that $\alpha_1$ is unique in $e$. Suppose that the color
$\alpha_1$ is assigned to at least two vertices $u,v \in e$, and
consider iteration $\alpha_1$ of $\A$. Let $H' = H[\{x \in V:
\chi(x) \geq \alpha_1\}]$, and let $\varphi$ be the $k$-colorful
coloring obtained for $H'$ in step $3$ of iteration $\alpha_1$.
Put $e' = \{ x \in e: \chi(x) \geq \alpha_1\}$. $e' \subset e$ is
a hyperedge in $H'$. Since $u,v \in e'$ then $|e| \geq 2$.
$\varphi$ is $k$-colorful for $H'$ so $e'$ contains at least two
vertices that are colored distinctively in $\varphi$. In
iteration $\alpha_1$, the vertices of one color class of
$\varphi$ are removed from $e'$. Since $e'$ contains vertices
from two color classes of $\varphi$, it follows that after
iteration $\alpha_1$ at least one vertex of $e'$ remains. Thus,
at least one vertex of $e'$ is colored in a later iteration than
$\alpha_1$, a contradiction to the maximality of $\alpha_1$. The
induction hypothesis is that in $\chi$ the colors
$\alpha_1,\ldots,\alpha_{j-1}$, $1 < j \leq k-1$,  all exist and
are unique in the hyperedge $e$. Consider the color $\alpha_j$.
There exists a vertex $u\in e$ such that $\chi(u)=\alpha_j$; for
otherwise it follows from the induction hypothesis that $|e| <
k-1$ since the colors $\alpha_1,\ldots,\alpha_{j-1}$ are all
unique in $e$ and $j-1 < k-1$. We prove that the color $\alpha_j$
is unique in $e$. Assume to the contrary that $\alpha_j$ is not
unique at $e$, and that in $\chi$ the color $\alpha_j$ is
assigned to at least two vertices $u,v \in e$. Put $H'' = H[\{u
\in V: \chi(u) \geq \alpha_j \}]$, and let $\varphi''$ be the
$k$-colorful coloring obtained for $H''$ in step $3$ of iteration
$\alpha_j$. Put $e''=\{u \in e: \chi(u) \geq \alpha_j\}$. $e''$
is a hyperedge of $H''$. By the induction hypothesis and the
definition of the colors $\alpha_1,\ldots,\alpha_{j-1}$, after
iteration $\alpha_j$ a set $U \subset e''$ of exactly $j-1$
vertices of $e''$ remains. In addition, $u,v \in e''$ and $U \cap
\{u,v\} = \emptyset$. Consequently, $|e''| \geq j+1$. Since
$\varphi''$ is $k$-colorful then $e''$ contains vertices from
$\min\{k,j+1\}$ color classes of $\varphi''$. $j \leq k-1$ so
$\min\{k,j+1\} = j+1$. Since in iteration $\alpha_j$ the vertices
of one color class of $\varphi''$ are removed from $e''$, it
follows that after iteration $\alpha_j$ at least $j$ vertices of
$e''$ remain. This is a contradiction to the induction
hypothesis.\QED

\noindent
\bfm{Remark.} Given a $k$-colorful coloring of $H$, the framework $\A$ obtains a
Strong Conflict-Free coloring of $H$ in a constructive manner. As mentioned above, in
this paper, computational efficiency is not of main interest.
However, it can be seen that for certain families of
geometrically induced hypergraphs, framework $\A$ produces an
efficient algorithm. In particular, for hypergraphs induced by
discs or axis-parallel rectangles, framework $\A$ produces an
algorithm with a low degree polynomial running time.
Colorful-colorings of such hypergraphs can be computed once the
arrangement of the discs is computed together with the depth of
every face (see, e.g., \cite{Sharir}). Due to space limitation we
omit the technical details involving the description of these
algorithms for computing $k$-colorful coloring for those
hypergraphs.

%\section{$\mbs{k}$-Strong-Conflict-Free Coloring of Geometrically Induced Hypergraphs}
\section{$\mbs{k}$-Strong-Conflict-Free Coloring of Geometrically Induced Hypergraphs}
\label{sec:bounds}

Theorems~\ref{th:reduce} and \ref{th:nonlinear} assert that in
order to attain upper bounds on $f_H(k)$, for a hypergraph $H$,
one may concentrate on attaining an upper bound on $c_H(k+1)$. In
this section we concentrate on colorful colorings.

\subsection{$\mbs{k}$-Strong-Conflict-Free Coloring and Union Complexity}
\label{sec:upper}

In this section, we prove Theorems~\ref{th:prior-res} and \ref{th:results}. Before proceeding with a proof of Theorem~\ref{th:prior-res}, we need several related definitions and theorems. A simple finite graph $G$ is called \emph{$k$-degenerate} if every vertex-induced sub-graph of
$G$ contains a vertex of degree at most $k$. For a finite set
$\R$ of simple closed planar Jordan regions and a fixed integer
$k$, let $G_k(\R)$ denote the graph with vertex set $\R$ and two
regions $r,s \in \R$ are adjacent in $G_k(\R)$ if there exists a
point $p \in \mbb{R}^2$ such that $(i)$ $p \in r \cap s$, and
$(ii)$ there exists at most $k$ regions in $\R \setminus \{r,s\}$
that contain $p$.

\begin{theorem}
\label{deg-prior} Let $\R$ be a finite set of simple closed
planar Jordan regions, let $H=H(\R)$, and let $k$ be a fixed
integer. If $G_k(\R)$ is $l$-degenerate then $c_H(k) \leq l+1 $.
\end{theorem}

Theorem~\ref{deg-prior} can be proved in a manner similar to that of Aloupis et al. (see
\cite{aloupis}) who proved Theorem~\ref{deg-prior} in the special case when $\R$ is a family of
discs. Due to space limitations, we omit a proof of this theorem. 

In light of Theorem~\ref{deg-prior}, in order to prove
Theorem~\ref{th:prior-res} it is sufficient to prove that for a
family of regions satisfying the conditions of
Theorem~\ref{th:prior-res} and a fixed integer $k$, the graph
$G_k(\R)$ is $O(k n^{\alpha})$-degenerate, where $\alpha$ is as in
Theorem~\ref{th:prior-res}.

\begin{lemma}
\label{degenerate} Let $k \geq 0$, let $0 \leq \alpha \leq 1$,
and let $c$ be a fixed constant. Let $\R$ be a set of $n$ simple
closed Jordan regions such that $\mcal{U}_{\R}(m) \leq c
m^{1+\alpha}$, for $1 \leq m \leq n$. Then $G_k(\R)$ is $O(k
n^{\alpha})$-degenerate.
\end{lemma}

Our approach to proving Lemma~\ref{degenerate} requires several
steps. These steps are described in the following lemmas. We shall
provide an upper bound on the average degree of every
vertex-induced subgraph of $G_k(\R)$ by providing an upper bound
on the number of its edges. We need the following lemma:

\begin{lemma}\bfm{(\cite{smoro})}
\label{smor} Let $\R$ be a set of $n$ simple closed planar Jordan
regions and let $\mcal{U}: \mbb{N} \rightarrow \mbb{N}$ be a
function such that $\mcal{U}(m)$ is the maximum union-complexity
of any $k$ regions in $\R$ over all $k \leq m$. Then the average
degree of $G_0(\R)$ is $O(\frac{\mcal{U}(n)}{n})$.
\end{lemma}

For a graph $G$, we write $E(G)$ to denote the set of edges of
$G$. We use Lemma~\ref{smor} to obtain the following easy lemma.
\begin{lemma}
\label{k=0} Let $0 \leq \alpha \leq 1$ and let $c$ be a fixed
constant. Let $\R$ be a set of $n$ simple closed Jordan regions
such that $\mcal{U}_{\R}(m) \leq c m^{1+\alpha}$, for $1 \leq m
\leq n$. Then there exists a constant $d$ such that $|E(G_0(\R))|
\leq \frac{dn^{1+\alpha}}{2}$.
\end{lemma}

\proof By Lemma~\ref{smor}, it follows that there exists a
constant $d'$ such that
$$
\frac{2 |E(G_0(\R))|}{n} =
\frac{ \sum_{x \in V(G_0(\R))} deg_{G_0(\R)}(x)}{n} \leq
\frac{d' c n^{1+\alpha}}{n}.
$$
Set $d = d'c$ and the claim follows.
\QED

For a set $\R$ of $n$ simple closed planar Jordan regions, define
$I(\R)$ to denote the graph whose vertex set is $\R$ and two
regions $r,s \in \R$ are adjacent if $r \cap s \not= \emptyset$.
The graph $I(\R)$ is called the \emph{intersection graph} of $\R$.
Note that for any integer $k$ $E(G_k(\R)) \subseteq E(I(\R))$.
Let $E \subseteq E(I(\R))$ be an arbitrary subset of the edges of
$I(\R)$. For every edge $e=(a,b) \in E$, pick a point $p_e \in
a\cap b$. Note that for distinct edges $e$ and $e'$ in $E$ it is
possible that $p_e = p_{e'}$. Put $X_{E,\R} = \{(p_e,r):\; e=(a,b) \in E \;
\mbox{and} \; r\in \R \setminus \{a,b\} \; \mbox{contains}\;
p_e\}$. In the following two lemmas we obtain a lower bound on
$|X_{E,\R}|$ in terms of $|E|$ and $|\R|$.

\begin{lemma}
\label{bootstrap} Let $0 \leq \alpha \leq 1$ and let $c$ and $d$
be the constants of Lemma~\ref{k=0}. Let $\R$ be a set of $n$
simple closed Jordan regions such that $\mcal{U}_{\R}(m) \leq c
m^{1+\alpha}$, for $1 \leq m \leq n$. Let $E \subseteq E(I(\R))$.
Then $|X_{E,\R}| \geq |E|-\frac{dn^{1+\alpha}}{2}$.
\end{lemma}

\proof Apply induction on the value
$|E|-\frac{dn^{1+\alpha}}{2}$. Let $P_E =\{p_e : e \in E\}$. One
may assume that $|E|-dn^{1+\alpha} \geq 0$ for otherwise the
claim follows trivially since $|X_{E,\R}| \geq 0$. Suppose
$|E|-\frac{dn^{1+\alpha}}{2} =1$. Since
$|E|>\frac{dn^{1+\alpha}}{2}$, then by Lemma~\ref{k=0} there
exists an edge $e=(a,b) \in E \setminus E(G_0(\R))$. Since $e
\notin E(G_0(\R))$, it follows that for every point $p \in a \cap
b$ there exists a region $r \in \R \setminus \{a,b\}$ such that
$p \in r$. Consequently, there exists a region $r \in \R
\setminus \{a,b\}$ such that $p_e \in r$. Hence, $(p_e,r) \in
X_{E,\R}$ and thus $|X_{E,\R}| \geq 1$. Assume that the claim
holds for $|E|-\frac{dn^{1+\alpha}}{2} =i$, where $i>1$, and
consider the case that $|E|-\frac{dn^{1+\alpha}}{2} =i+1$. Let $e
=(a,b) \in E$ be an edge such that there exists a region $r \in \R
\setminus \{a,b\}$ with $p_e \in r$. Define $E' = E \setminus
\{e\}$. Note that $P_{E'} \subset P_E$ and
$|E'|-\frac{dn^{1+\alpha}}{2} =i$. By the induction hypothesis it
follows that $|X_{E',\R}| \geq |E'|-\frac{dn^{1+\alpha}}{2}$.
Observe that $X_{E',\R} \subset X_{E,\R}$ and that $|X_{E,\R}|
\geq |X_{E',\R}|+1$. It follows that
$$|X_{E,\R}| \geq |X_{E',\R}|+1 \geq |E'|-\frac{dn^{1+\alpha}}{2} +1 = i+1 = |E|-\frac{dn^{1+\alpha}}{2}.$$ \QED

\begin{observation}
\label{binom} Let $0 \leq \alpha \leq 1$ and let $X$ be a
binomial random variable with parameters $n$ and $p$. Then
$$
\expect{X^{1+\alpha}} \leq \expect{Xn^{\alpha}} =
n^{\alpha}\expect{X} =n^{1+\alpha}p.
$$
\end{observation}

\begin{lemma}
\label{prob} Let $0 \leq \alpha \leq 1$ and let $c$ and $d$ be
the constants of Lemma~\ref{k=0}. Let $\R$ be a set of $n$ simple
closed Jordan regions such that $\mcal{U}_{\R}(m) \leq c
m^{1+\alpha}$ for $1 \leq m \leq n$. Let $E \subseteq E(I(\R))$
such that $|E| > d n^{1+\alpha}$ and let $\{p_e | e \in E\}$ and
$X_{E,\R}$ be as before. Then $|X_{E,\R}| \geq \frac{|E|^2}{2 d
n^{1+\alpha}}$.
\end{lemma}

\proof  Let $\R' \subseteq \R$ be a subset of regions of $\R$
chosen randomly and independently such that for every region $r
\in \R$, $\prob{r \in \R'} = p$ for $p =
\frac{dn^{1+\alpha}}{|E|}$ (note that $p < 1$). Let $E' \subseteq
E$ be the subset of edges that is defined by the intersections of
regions in $\R'$. Let $P_{E'} = \{p_e: e \in E'\}$. $P_{E'}
\subseteq P_E$ and thus $X_{E',\R'} \subseteq X_{E,\R}$. Each of
$|\R'|,|E'|$, and $|X_{E',\R'}|$ is a random variable.

By Lemma~\ref{bootstrap} and by linearity of expectation, it
follows that $\expect{|X_{E',\R'}|} \geq \expect{|E'|} -
\expect{\frac{d}{2} |\R'|^{1+\alpha}}$.
% = \expect{|E'|} - \frac{d}{2} \expect{\left( \sum_{i=1}^n X_i \right)^{1+\alpha}}$.

By Observation~\ref{binom}, $\expect{\frac{d}{2}
|\R'|^{1+\alpha}} \leq \frac{d}{2}n^{1+\alpha}p$. Hence, it
follows that $\expect{|X_{E',\R'}|} \geq \expect{|E'|} -
\frac{d}{2}n^{1+\alpha}p$. For an edge $e = (a,b) \in E$,
$\prob{e \in E'} = \prob{ a,b \in \R' } = p^2$ so $\expect{|E'|} =
p^2 |E|$. In addition, for an edge $e =(a,b)$ and a region $r \in
\R \setminus \{a,b\}$ $\prob{(p_e,r) \in X_{E',\R'}} =
\prob{a,b,r \in \R'} = p^3$. Thus, $\expect{|X_{E',\R'}|} = p^3
|X_{E,\R}|$. It follows that $|X_{E,\R}| \geq \frac{|E|}{p} -
\frac{dn^{1+\alpha}}{2 p^2}$. Substituting the value of $p$ in
the latter inequality completes the proof of the lemma.\QED

Next, a proof of Lemma~\ref{degenerate} is presented.\\

\noindent
\bfm{Proof of Lemma~\ref{degenerate}.} Let $d$ be the constant
from Lemma~\ref{k=0}. Let $V \subseteq V(G_k(\R))$ be a subset of
of $m$ vertices and let $G$ be the subgraph of $G_k(\R)$ induced
by $V$. Define $E = E(G)$. Observe that $E \subseteq E(I(\R))$.
There are two cases: Either $|E| \leq d m^{1+\alpha}$ or $|E| > d
m^{1+\alpha}$. In the former case, the average degree of a vertex
in $G$ is at most $2dm^{\alpha}$. In the latter case, it follows
from Lemma~\ref{prob} that $|X_{E,V}| \geq
\frac{|E|^2}{2dm^{1+\alpha}}$. On the other hand, since $E
\subseteq E(G_k(\R))$ then by definition, for every edge $e \in E$
the chosen point $p_e$ can belong to at most $k$ other regions of
$\R$. Thus $|X_{E,V}| \leq k|E|$. Combining these two
inequalities we have: $|E| \leq 2dk m^{1+\alpha}$ and thus the
average degree of $G$ in this case is at most $4dkm^{\alpha}$.
Hence, in $G$ there exists a vertex whose degree is at most $\max
\{2dm^{\alpha},4dkm^{\alpha} \} = 4dkm^{\alpha}$. \QED

As mentioned in the introduction, Theorem~\ref{th:results} is a corollary of a combination of
Theorems~\ref{th:reduce}, \ref{th:nonlinear}, and Theorem~\ref{th:prior-res}.

\subsection{$\mbs{k}$-Strong Conflict-Free Coloring of Axis-Parallel Rectangles}
\label{sec:rects}

In this section, we consider $kSCF$-colorings of axis-parallel
rectangles and prove Theorems~\ref{th:rects} and \ref{th:colorful_rects}. 
As mentioned in the introduction, a proof of Theorem~\ref{th:rects} can be derived from a combination of Theorem~\ref{th:reduce} and Theorem~\ref{th:colorful_rects}. Consequently, we concentrate on a proof of Theorem~\ref{th:colorful_rects}. To that end we require the following lemma. 

\begin{lemma}
\label{lem:line} Let $k \geq 2$. Let $\R$ be a set of $n$
axis-parallel rectangles such that all rectangles in $\R$
intersect a common vertical line $\ell$, and let $H=H(\R)$. Then
$c_H(k) =O(k)$.
\end{lemma}

\proof Assume, without loss of generality, that the rectangles are
in general position (that is, no three rectangles' boundaries
intersect at a common point). According to
Theorem~\ref{th:prior-res}, it is sufficient to prove that for
every subset of rectangles $\R' \subseteq \R$, the
union-complexity of $\R'$ is at most $O(|\R'|)$. Let $\R'
\subseteq \R$ and consider the boundary of the union of the
rectangles of $\R'$ that is to the right of the line $\ell$. Let
$\partial_r \R'$ denote this boundary. An intersection point on
$\partial_r \R'$ results from the intersection of a horizontal
side of a rectangle and a vertical side of another rectangle.
Each horizontal side of a rectangle in $\R'$ may contribute at
most one intersection point to $\partial_r \R'$. Indeed, let $s$
be a horizontal rectangle side. Let $p$ be the right-most
intersection point on $s$ to the right of the line $\ell$ and let
$q$ be any other intersection point on $s$ to the right of
$\ell$. Let $r$ be the rectangle whose vertical side defines $p$
on $s$. Since $r$ intersects $\ell$, the point $p$ lies on the
right vertical side of $r$. Hence, $q \in r$; for otherwise
either $r$ does not intersect $\ell$ or $q$ is to the left of
$\ell$, in which case $q$ does not lie on $\partial_r \R'$. It
follows that every horizontal side $s$ of some rectangle
contributes at most one point to $\partial_r \R'$. As there are
$2|\R'|$ such sides then $\partial_r \R'$ contains at most
$2|\R'|$ points. A symmetric argument holds for the boundary of
$\partial \R'$ that lies to the left of $\ell$. Hence, the
union-complexity of $\R'$ is at most $O(|\R'|)$. By
Theorem~\ref{th:prior-res}, the claim follows. \QED

Next, we prove Theorem~\ref{th:colorful_rects}.\\

\noindent
\bfm{Proof of Theorem~\ref{th:colorful_rects}.} Let $\ell$ be a vertical line such that at most $n/2$ rectangles lie fully to its right and to its left, respectively.
Let $\R'$ and $\R''$ be the sets of rectangles that lie entirely
to the right and entirely to the left of $\ell$, respectively.
Let $\R_{\ell}$ denote the set of rectangles in $\R$ that
intersect $\ell$, and let $c(n)$ denote the least number of
colors required by a colorful coloring of any $n$ axis-parallel
rectangles. By Lemma~\ref{lem:line}, the set of rectangles
$\R_{\ell}$ can be colored using $O(k)$ colors. In order to
obtain a $k$-colorful coloring of $\R$, we color $\R_{\ell}$ using
a set $D$ of $O(k)$ colors. We then color $\R'$ and $\R''$
recursively by using the same set of colors $D'$ such that $D\cap
D' = \emptyset$. The function $c(n)$ satisfies the recurrence
$c(n) \leq O(k) +c(n/2)$. Thus, $c(n) = O(k \log n)$. Let
$\varphi$ be the resulting coloring of the above coloring
procedure. It remains to prove that $\varphi$ is a valid
$k$-colorful coloring of $\R$. The proof is by induction on the
cardinality of $\R$. Suppose $\R'$ and $\R''$ are colored
correctly under $\varphi$, and consider a point $p \in \bigcup_{r
\in \R} r$. If $r(p) \subset \R_{\ell}$ or $r(p) \subset \R'$ or
$r(p) \subset \R''$ then by Lemma~\ref{lem:line} and the induction
hypothesis, $r(p)$ is colored correctly under $\varphi$. It is
not possible that both $r(p) \cap \R' \not= \emptyset$ and $r(p)
\cap \R'' \not= \emptyset$. Hence, it remains to consider points
$p$ for which either $r(p) \subset \R_{\ell} \cup \R'$ or $r(p)
\subset \R_{\ell} \cup \R''$. Consider a point $p$ which is,
w.l.o.g, of the former type. Let $i=|r(p) \cap \R_{\ell}|$ and
$j=|r(p) \cap \R'|$. If either $i\geq k$ or $j \geq k$, then
either by Lemma~\ref{lem:line} or by the inductive hypothesis the
hyperedge $r(p)$ is $k$-colorful. It remains to consider the case
that $i+j \geq k$ and $i,j <k$. Let $\varphi_{\ell}$ and
$\varphi_{\R'}$ be the colorings of $\R_{\ell}$ and $\R'$ induced
by $\varphi$, respectively. By the inductive hypothesis, the
rectangles in the set $r(p) \cap \R'$ are colored distinctively
using $j$ colors under $\varphi_{\R'}$. In addition, by
Lemma~\ref{lem:line}, the rectangles in $r(p) \cap \R_{\ell}$ are
colored using $i$ distinct colors under $\varphi_{\ell}$.
Moreover, the colors used in $\varphi_{\ell}$ are distinct from
the ones used in $\varphi_{\R'}$. Hence, $r(p)$ is $\min
\{|r(p)|,k\}$-colorful. This completes the proof of the lemma.
\QED

\bibliographystyle{abbrv}
\bibliography{refs}

\end{document}